\documentclass[12pt]{article}
\usepackage[latin1]{inputenc}
\usepackage{amssymb,amsmath,amsfonts}
\newtheorem{theorem}{Theorem}
\newtheorem{proposition}[theorem]{Proposition}
\newtheorem{lemma}[theorem]{Lemma}

\newcommand{\R}{\mathbb{R}}

\newcommand{\Sf}{\mathbb{S}}
\newcommand{\C}{\mathbb{C}}

\newcommand{\spa}{\mbox{span}}

\newcommand{\po}{{\hspace*{-1ex}}{\bf .  }}

\newcommand{\nap}{\nabla^{\perp}}

\newcommand{\nab}{\tilde\nabla}

\def\<{{\langle}}
\def\>{{\rangle}}

\def\n{\nabla}
\def\d{\partial}

\def\a{\alpha}
\def\va{\varphi}

\def\bea{\begin{eqnarray*} }
\def\eea{\end{eqnarray*} }
\def\be{\begin{equation} }
\def\ee{\end{equation} }

\def\nap{\nabla^\perp}
\def\proof{\noindent\emph{Proof: }}
\def\qed{\ifhmode\unskip\nobreak\fi\ifmmode\ifinner
\else\hskip5 pt \fi\fi\hbox{\hskip5 pt \vrule width4 pt
height6 pt  depth1.5 pt \hskip 1pt }}

\begin{document}

\title{A class of minimal submanifolds in spheres}
\author{M.\ Dajczer and Th.\ Vlachos}

\date{}
\maketitle
\renewcommand{\thefootnote}{\fnsymbol{footnote}} 
\renewcommand{\thefootnote}{\arabic{footnote}} 

\renewcommand{\thefootnote}{\fnsymbol{footnote}} 
\renewcommand{\thefootnote}{\arabic{footnote}} 

\begin{abstract} We introduce a class of minimal submanifolds $M^n$, $n\geq 3$, 
in spheres $\Sf^{n+2}$ that are ruled by totally geodesic spheres of dimension $n-2$.  
If simply-connected, such a submanifold admits a one-parameter associated family 
of equally ruled minimal isometric deformations that are genuine. 
As for compact examples, there are plenty of them but only for dimensions $n=3$ 
and $n=4$.  In the first case, we have that $M^3$ must be a  
$\Sf^1$-bundle over a minimal torus $T^2$ in $\Sf^5$ and in the second case  
$M^4$ has to be a $\Sf^2$-bundle over a minimal sphere $\Sf^2$ in $\Sf^6$.  
In addition, we provide new examples in relation to the well-known 
Chern-do Carmo-Kobayashi problem since taking the torus $T^2$ to be flat yields
a minimal submanifolds $M^3$ in $\Sf^5$ with constant scalar curvature.
\end{abstract} 

In several directions, this paper should be considered as a continuation of our 
work in \cite{dv1} where a new class of minimal ruled submanifolds $M^n$ of
Euclidean space $\R^{n+2}$, $n\geq 3$, were studied.  
These submanifolds lay in codimension two and may be metrically 
complete regardless the dimension. The rulings are of codimension 
two in the manifold whereas the rank, that is, the complement 
of the index of relative nullity, is $\rho=4$  (unless $n=3=\rho$) along 
an open dense subset.
If simply-connected, the submanifolds admit a $\Sf^1$-parameter family of 
genuine isometric deformations. 
Hence, this class of examples should be seen as a new addition to the possible, 
local or global, classification of Euclidean submanifolds in codimension 
two that admit genuine isometric deformations; see \cite{dv1} for a
discussion of that open problem.

In this paper, we consider a similar construction but for the round sphere as 
ambient space. We obtain minimal submanifolds $M^n$ in $\Sf^{n+2}$, $n\geq 3$,
with similar properties as the ones in the Euclidean space.
Notice that being ruled now 
means that the submanifold carries a foliation by (open subsets of) totally 
geodesic spheres in $\Sf^{n+2}$ of dimension $n-2$. 
If the manifold is simply-connected, by taking the cones in $\R^{n+3}$ of the 
components  in the associated family in $\Sf^{n+2}$ we obtain a new class of 
genuinely deformable Euclidean submanifolds in codimension two but, of course, 
these are not complete.

New examples of minimal submanifolds in spheres are certainly welcome since the 
explicitly known ones are usually quite elaborate and certainly less abundant than 
in the Euclidean space.  Frequently, they are spheres of constant sectional curvature
or products of them. On the other hand, the submanifolds here introduced can 
only be complete, or even compact, for dimensions $n=3$ or $4$. 
If compact and according to the dimension, the submanifold must be topologically 
either a $\Sf^1$-bundle over $\Sf^1\times\Sf^1$ in $\Sf^5$ or a $\Sf^2$-bundle 
over $\Sf^2$ in $\Sf^6$.

The compact examples in the case of the torus are of particular interest 
by two quite different reasons. First, if we replace the torus by its universal cover
we obtain a three-dimensional manifold  that is not longer compact but has
an $\Sf^1$-parameter family of isometric minimal deformations. But  
the compact submanifold itself only admits, at most, a finite set of isometric minimal 
deformations. The second reason, has to do with the well-known 
Chern-do Carmo-Kobayashi problem \cite{CCK} that concerns compact 
minimal submanifolds of the sphere with constant scalar curvature.
We show that if the torus considered is flat, and these were all parametrically 
described by Miyaoka \cite{M}, then  $M^3$ has constant scalar curvature.

\section{The results}

This section is devoted to state the results of the paper while proofs 
are left for the following one.  Up to the last two results, the other 
theorems in this paper can be seen as the ``spherical version" of the 
results obtained in \cite{dv1} for submanifolds in Euclidean space
\vspace{1,5ex}

Let $g\colon L^2\to\Sf^{n+2}$ denote a substantial oriented  minimal surface. 
As already recalled in \cite{dv1} the normal 
bundle $N_gL$ of $g$ splits along an open dense subset of $L^2$ as
$$
N_gL=N_1^g\oplus N_2^g\oplus\dots\oplus N_m^g,\;\;\; m=[(n+1)/2],
$$
where each subbundle $N_s^g$, $1\leq s\leq m$,  is spanned by the corresponding 
$(s+1)^{th}$-fundamental form $\a_g^{s+1}\colon TL\times\cdots\times TL\to N_gL$
and has rank two except possible the last one that has rank one if $n$ is odd.

If  $L^2$ is simply-connected,  there exists a one-parameter associated family 
of minimal isometric immersions. In fact, for each constant 
$\theta\in\Sf^1=[0,\pi)$ consider the parallel orthogonal tensor field 
$$
J_{\theta}=\cos\theta I+\sin\theta J
$$
where $I$ is the identity map and $J$ the complex structure determined by the
metric and orientation.  Then, the symmetric section $\a_g(J_\theta\cdot,\cdot)$ 
of the bundle $\text{Hom}(TL\times TL,N_g L)$ satisfies the Gauss, Codazzi and Ricci 
equations with respect to the same induced normal connection; 
see \cite{dg} for details. Therefore, there exists an isometric minimal immersion  
$g_{\theta}\colon L^2\to\Sf^{n+2}$  whose second fundamental form is
$$
\a_{g_{\theta}}(X,Y)=\phi_\theta\a_g(J_{\theta}X,Y)
$$ 
where $\phi_\theta\colon N_gL\to N_{g_{\theta}}L$ is the parallel 
vector bundle isometry that identifies the normal bundles  as well as
each normal subbundles $N_s^g$ with $N_s^{g_\theta}$ for any $1\leq s\leq m$.
\vspace{1ex}

In the sequel, let $g\colon L^2\to\Sf^{n+2}$, $n\geq 2$ be a substantial  
\emph{$1$-isotropic} surface.  
This means that $g$ is minimal and that the  
ellipse of curvature (of first order) at any point is a circle. 
Let $L_0$ be the open subset of $L^2$ where $\dim N_{1}^g(p)=2$. 
It was shown in \cite{dv} that $L^2\smallsetminus L_0$ consists of 
isolated points and that the vector  bundle $N_1^g|_{L_0}$  smoothly extends 
to a plane bundle over $L^2$, that we still denote by $N^g_1$. 
 
Let $\pi\colon\Lambda_g\to L^2$ denote the vector 
bundle of rank $n-2$ whose fibers are the orthogonal complement in the normal 
bundle $N_gL$ of $g$ of its extended first normal bundle $N_1^g$.  
Then $F_g\colon M^n\to\Sf^{n+2}$ is the submanifold of $\Sf^{n+2}$ 
associated to $g$ constructed by attaching at each point of the surface $g$ 
the totally geodesic  sphere $\Sf^{n-2}$ whose tangent space at that point 
is the  fiber  of $\Lambda_g$, that is, 
\be\label{param}
(p,v)\in\Lambda_g\mapsto F_g(p,v)=\exp_{g(p)}v,
\ee
while dropping the singular points whenever they exist, i.e., 
points where the induced metric is singular.
By definition $F_g$ is an $(n-2)$-ruled submanifold, that is, there is an 
integrable tangent distribution of dimension $n-2$ whose leaves are mapped 
diffeomorphically by $F_g$ onto  open subsets of totally geodesic  
$(n-2)$-spheres of $\Sf^{n+2}$.

For simplicity, it is very convenient to do computations in terms of 
the cone of $M^n$ in $\Sf^{n+2}\subset\R^{n+3}$, and then view $M^n$  
as the intersection of that cone with $\Sf^{n+2}$. More precisely, 
we consider the map $G_g\colon\R\times\Lambda_g\to\R^{n+3}$ given by
\be\label{map}
G_g(s,p,v)=sg(p)+v
\ee
and set 
$S_G=\{0\}\times\big(\Lambda_g^*(p)\smallsetminus\{0\}\big)$
where
$$
\Lambda_g^*(p)=\{(p,v)\in\Lambda_g:v\perp N_2^g(p)\}.
$$
In next section, we show that the set of singular points of the metric induced
by $G_g$ consists of the vertex  $V=(0,p,0)$  and the set $S_G$.  Set 
$$
N^{n+1}=\R\times\Lambda_g\smallsetminus (V\cup S_G)
$$
and denote $G=G_g|_{N^{n+1}}$. Thus, we have
$$
M^n=\{(s,p,v)\in\R\times\Lambda_g\smallsetminus S_G: s^2+\|v\|^2=1\}
$$
and $F_g=G|_{M^n}$ where $M^n$ is endowed with the induced metric.
Observe that $M^n$ is complete (respectively, compact) if and only if $g$ is 
complete (respectively, compact) and $S_G$ is empty.
Notice also that $S_G$ can only be empty for $n=3,4$.
\vspace{1ex}

In the sequel, we denote by ${\cal H}$ the tangent distribution 
orthogonal to the rulings.
An embedded surface $j\colon L^2\to M^n$ is called an \emph{integral surface} 
of ${\cal H}$ if $j_*T_pL={\cal H}(j(p))$ at every point $p\in L^2$.

\begin{theorem}\po\label{main1} Let $g\colon L^2\to\Sf^{n+2}$, $n\geq 3$, be a 
$1$-isotropic substantial surface. Then the associated immersion 
$F_g\colon M^n\to\Sf^{n+2}$ is an $(n-2)$-ruled minimal 
submanifold with rank  $\rho=4$ (unless $n=3=\rho$) on an open dense 
subset of $M^n$. Moreover,  the integral surface $L^2$ of ${\cal H}$ is totally 
geodesic and unique up to the one obtained by composing with the 
antipodal map. 

Conversely, let  $F\colon M^n\to\Sf^{n+2}$ be an $(n-2)$-ruled 
minimal immersion with $n\geq 4$ and $\rho=4$ (unless $n=3=\rho$) 
on an open dense subset of $M^n$.
Assume that ${\cal H}$ admits a  totally geodesic integral surface 
$j\colon L^2\to M^n$ which is a global cross section to the rulings. 
Then the surface $g=F\circ j\colon L^2\to\Sf^{n+2}$ is $1$-isotropic and 
$F$ can be parametrized as $F_g$.
\end{theorem}

The existence of genuine deformations is considered in the following result.

\begin{theorem}\po\label{main2} Let $g\colon L^2\to\Sf^{n+2}, n\geq 3$, 
be a simply-connected $1$-isotropic substantial surface. Then $F_g$ allows 
a smooth one-parameter family of minimal genuine isometric deformations
$F_\theta\colon M^n\to\Sf^{n+2},\;\theta\in\Sf^1$, 
such that $F_0=F_g$ and each $F_\theta$ carries the same rulings and 
relative nullity leaves as $F_g$. 
\end{theorem}

The relation between the second fundamental forms of members
of the associated family is given next, for simplicity, in terms of their cones.

\begin{theorem}\po\label{af}  Let $g\colon L^2\to\Sf^{n+2}, n\geq 3$, 
be a simply-connected $1$-isotropic substantial surface. Then $G$ allows 
an associated smooth one-parameter family of minimal genuine isometric immersions 
$G_\theta\colon N^{n+1}\to\R^{n+3},\;\theta\in \Sf^1$, such that $G_0=G$ 
and each $G_\theta$ carries the same rulings and relative nullity leaves 
as $G$.

Moreover, there is  a parallel vector bundle isometry 
$\Psi_\theta\colon N_GN\to N_{G_\theta}N$
such that the relation between the second 
fundamental forms  is given by
\be\label{forms}
\a_{G_\theta}(X,Y)=\Psi_\theta \big(\textsf{R}_{-\theta} \a_G(X,Y) 
+2\kappa\sin(\theta/2)\beta ({\cal J}_{-\theta /2}X,Y)\big)
\ee
where $\textsf{R}_\theta$ is the rotation of angle
$\theta$ on $N_GN$ that preserves orientation,
$\kappa$ is the radius of the ellipse of curvature of $g$
and  $\beta$ is the traceless bilinear form  defined by (\ref{b}).
\end{theorem}

A substantial surface in even codimension $g\colon L^2\to\Sf^{n+2}$ 
is called \emph{pseudoholomorphic} when the ellipses of curvature of 
any order are circles at any point.  
In odd codimension, the surface is called \emph{isotropic} 
when the ellipses of curvature of any order but for the
last one-dimensional normal subbundle are circles 
at any point. 
\vspace{1,5ex}

If $g\colon L^2\to\Sf^{n+2}$ is pseudoholomorphic, then taking a rotation of 
angle $\theta\in\Sf^1$ 
that  preserves orientation in each $N_s^g$, $s\geq 2$, induces an intrinsic 
isometry $S_\theta$ on $M^n$. The next result says that $F_g$ is equivariant 
with respect to the one-parameter family of intrinsic isometries $S_\theta$.

\begin{theorem}\po\label{main3}  If $g\colon L^2\to\Sf^{n+2}$ is 
pseudoholomorphic, then $F_g\circ S_{-\theta}$ is congruent to 
$F_{\theta}$ for any $\theta\in\Sf^1$.  
\end{theorem}

We have that $F_g\colon M^3\to\Sf^5$ or $F_g\colon M^4\to\Sf^6$ is compact
if and only if $L^2$ is compact and $g$ is regular. The latter condition 
means that $L_0$ is empty and that $N_2^g$ has constant dimension. 
According to a result of Asperti \cite{A} any 
compact regular substantial minimal surface in $\Sf^5$ is a topological 
torus and in $\Sf^6$ is a topological sphere. For both cases, there are plenty 
of $1$-isotropic examples. 
In fact, the  tori in $\Sf^5$ include the flat ones described 
parametrically by Miyaoka \cite{M} and those that are holomorphic 
with respect to the nearly Kaehler structure of $\Sf^6$ considered in  
\cite{Br}, \cite{EV} and \cite{H1}. Other examples of $1$-isotropic
surfaces in $\Sf^5$ are the Legendrian surfaces given in {\cite{U}}.

Minimal $2$-spheres in spheres have been investigated 
by Calabi, Barbosa and Chern among others.  From their work, we know that
these surfaces must be substantial in even codimension and pseudoholomorphic.  
It was then shown by Calabi \cite{Cal} that 
any such surface in $\Sf^6$ is regular if its area is $24\pi$. 
Then Barbosa \cite{B} proved that the space of  these 
surfaces is diffeomorphic to  $SO(7,\C)/SO(7,\R)$, where  $SO(7,\C)$
denotes the set of $7\times 7$ complex matrices that satisfy $AA^t = I$ 
and $\det A=1$. \vspace{1,5ex}

Concerning the set of genuine minimal isometric deformations of compact 
submanifolds constructed from tori we have the following result.

\begin{theorem}\po\label{mod} Let $g\colon L^2\to\Sf^5$ be a regular substantial 
isotropic surface. Then, the set of all equally ruled minimal isometric immersions 
of $M^3$ into $\Sf^5$ as $F_g\colon M^3\to\Sf^5$  is finite or parametrized by a circle $\Sf^1$. 
If $L^2$ is compact then the set is necessarily finite.
\end{theorem}

As discussed in the introduction the last result is of independent interest.
 
\begin{theorem}\po\label{last} Let $g$  be a flat $1$-isotropic torus in $\Sf^5$. 
Then $F_g\colon M^3\to\Sf^5$ is a compact minimal submanifold with constant normalized scalar curvature 
$s=-1/3$.
\end{theorem}

\section{The proofs}

In this section, we provide several proofs for $n\geq 4$ but similar 
arguments take care of the case $n=3$.\vspace{1,5ex}

First we discussed the set of singular points of $F_g$.

\begin{proposition}\po\label{sp} Let  $g\colon L^2\to\Sf^{n+2},n\geq 4$, 
be a substantial oriented minimal surface. Then, the set of singular points 
of the map $G\colon\R\times\Lambda_g\to\R^{n+3}$ given by 
(\ref{map}) consists of $V=(0,p,0)$ and $S_G$.
\end{proposition}

\proof
Fix $(s_0,p_0,v_0)\in\R\times\Lambda_g \smallsetminus\{V\}$. Choose a smooth 
orthonormal frame $\{e_5,\dots,e_{n+2}\}$ of $\Lambda_g$ on a neighborhood $U$ 
of $p_0$ and set
$$
v_0=\sum_{i\geq1}a_{i}e_{i+4}(p_0).
$$
Consider the  projection $\Pi\colon \R\times\Lambda_g \to L^2$ and parametrize
$\Pi^{-1}(U)$
via the diffeomorphism $h\colon U\times\R^{n-1}\to\Pi^{-1}(U)$ given by
$$
h(p,s,t_1,\dots, t_{n-2})=\big(s,p,\sum_{i\geq1}t_{i}e_{i+4}\big).
$$
That $(s_0,p_0,v_0)\in S_G$ means that there exists a non-zero vector
$$
Z=X+\lambda_0 \partial/\partial s
+\sum_{i\geq1}\lambda_{i}\partial/\partial t_i 
\in\ker (G\circ h)_*(p_0,s_0,a_1,\dots,a_{n-2})
$$
where $X\in T_{p_0}L$. Thus, 
$$
\lambda_0 g(p_0)+s_0g_*(p_0)X+\sum_{i\geq1}a_i\nap_Xe_{i+4}(p_0)
+\sum_{i\geq1}\lambda_{i}e_{i+4}(p_0)=0.
$$
Since $Z\neq 0$, we obtain that $\lambda_0=0$, $s_0=0$, $X\neq 0$ and
$$
\sum_{i\geq1}a_i\nap_Xe_{i+4}(p_0)+\sum_{i\geq1}\lambda_{i}e_{i+4}(p_0)=0.
$$
It follows that
$$
\<v_0,\nap_X\xi\>(p_0)=0
$$
for any $\xi\in N_1^g$. We easily conclude that $v_0\perp N_2^g(p_0)$. 
The converse is immediate. 
\vspace{1,5ex}\qed

In the sequel, we argue for an open set of $L^2$ where all the normal
subspaces  $N^g_s$'s of the substantial 
oriented minimal surface $g\colon L^2\to\Sf^{n+2}$ have constant dimension.
Choose local positively oriented orthonormal frames $\{e_1,e_2\}$ in $TL$ and 
$\{e_3,e_4\}$ of $N_1^g$ such that 
$$
\a_g(e_1,e_1)=\kappa e_3\;\;\;\mbox{and}\;\;\; \a_g(e_1,e_2)=\mu e_4
$$ 
where $\kappa,\mu$ are the semi-axes of the ellipse of curvature. Take 
a local orthonormal normal frame $\{e_5,\ldots,e_{n+2}\}$ 
such that $\{e_{2r+1},e_{2r+2} \}$ is positively oriented spanning 
$N_r ^g$ for every even $r$. When $n=2m+1$ is odd, then $e_{2m+1}$ spans the last 
normal bundle. We refer to  $\{e_1,\ldots,e_{n+2}\}$ as an \emph{adapted frame} 
of $g$ and consider the one-forms
$$
\omega_{ij}=\<\nab e_i,e_j\>\;\;\mbox{for}\;\;1\leq i,j\leq n+2,
$$
where $\nab$ denotes the Riemannian connection in the ambient space. 
Using that
$$
\a_g^3(e_1,e_1,e_1)+\a_g^3(e_1,e_2,e_2)=0
$$
we easily obtain
\be\label{conn}
\omega_{45}=-\dfrac{1}{\lambda}*\omega_{35}
\;\;\mbox{and}\;\;\omega_{46}=-\dfrac{1}{\lambda}*\omega_{36}
\ee
where $\lambda=\mu/\kappa$ and $*$ denotes the Hodge operator, i.e., 
$*\omega(e)=-\omega(Je)$. Here  $J$ is the  complex structure of $L^2$ induced
by the orientation.  We denote by 
$$
V=a_1e_1+a_2e_2,\;\; W=b_1e_1+b_2e_2,\;\;
Y=c_1e_1+c_2e_2\;\;\mbox{and}\;\;
Z=d_1e_1+d_2e_2
$$
the dual vector fields of $\omega_{35},\omega_{36},\omega_{45}$ and $\omega_{46}$, 
respectively. Then (\ref{conn}) is equivalent to 
$$
Y=-\frac{1}{\lambda}JV \;\;\mbox{and}\;\; Z=-\frac{1}{\lambda}JW,
$$
and hence
\be\label{omegas}
\lambda c_1 = a_2,\;\; \lambda c_2=-a_1,\;\; 
\lambda d_1 = b_2\;\;\mbox{and}\;\;\lambda d_2=-b_1.
\ee

Clearly, we have that $G\colon N^{n+1}\to\R^{n+3}$  is an immersion and 
$$
T_{(s,p,v)}N=\R\oplus T_{(p,v)}\Lambda_g
=\R\oplus\mathcal H^G(p,v)\oplus\mathcal V(p,v)
$$ 
where $\R=\spa\{\partial/\partial s\}$ and $\mathcal H^G$ is the 
orthogonal complement of $\mathcal V$ in $T\Lambda_g$. Moreover, 
$\mathcal V$ denotes the vertical bundle of $\pi\colon\Lambda_g\to L^2$ 
given by $\mathcal V=\ker\pi_*$.

Fixed  $(p,v)\in\Lambda_g$, let $\delta_v$ be the normal vector field
defined in a neighborhood of $p$ by
\be\label{delta}
\delta_v(q)=\sum_{j\geq 5}\<v,e_j(p)\> e_j(q).
\ee 
Let $\beta_i$, $1\leq i\leq 2$, be the curves in $\Lambda_g$ satisfying
$\beta_i(0)=(p,v)$ given by 
$$
\beta_i(t)=(c_i(t),\delta_v(c_i(t))) 
$$
where $c_i(t)$ is a smooth curve in a neighborhood of $p$ satisfying 
$c_i'(0)=e_i(p)$. Set 
\be\label{yes}
Y_i=\beta_i'(0)\in T_{(p,v)}\Lambda_g,\;1\leq i\leq 2.
\ee
Let $G_i,H_i\in C^\infty(\Lambda_g)$, $\;1\leq i\leq 2$, be the functions 
$$
G_i=t_2\omega^i_{56}+t_3\omega^i_{57}+t_4\omega^i_{58},\;\;
H_i=-t_1\omega^i_{56}+t_3\omega^i_{67}+t_4\omega^i_{68}
$$
where $\omega_{ij}^k=\omega_{ij}(e_k)$ and 
$t_j\in C^\infty(\Lambda_g)$ is defined by
$$
t_j(q,w)=\<w,e_{j+4}(q)\>,\;\; 1\leq j\leq 4.
$$

It is clear that $G_*(s,p,v)\mathcal V=(N_1^g(p))^\perp\subset N_gL(p)$ 
holds up to parallel identification in $\R^{n+3}$.
The vector bundle $\mathcal V$ can be orthogonally decomposed as
$\mathcal{V}=\mathcal{V}^1\oplus \mathcal{V}^0$
where $\mathcal{V}^1$ denotes the plane bundle determined by 
$$
G_*(s,p,v)\mathcal V^1=N_2^g (p).
$$
Let $\{E_3,E_4\}$ and $\{E_5,\dots, E_n\}$ be local orthonormal frames  
of $\mathcal V^1$ and  $\mathcal V^0$, respectively, such that 
$$
G_* E_j= e_{j+2}\;\;\mbox{for}\;\; 3\leq j\leq n.
$$

\begin{lemma}\label{NF}\po The vectors  $X_1,X_2\in T_{(p,v)}\Lambda_g$ 
defined as  
\be\label{X}
X_i= Y_i+ G_iE_3+H_iE_4-\sum_{j\geq 7}\<\nap_{e_i}\delta_v,e_j\>E_{j-2}
\ee
satisfy that $X_1,X_2\in\mathcal H^G(p,v)$ and that
$$
G_*X_1=sg_* e_1 -\va_1e_3
-\dfrac{1}{\lambda}\va_2e_4,\;\;\;
G_*X_2=sg_* e_2 -\va_2e_3
+\dfrac{1}{\lambda}\va_1e_4
$$
where $\va_j=t_1^0a_j+t_2^0b_j$ and $t_j^0=t_j(p,v)$. Moreover, 
the space $N_GN(s,p,v)$ is spanned by 
$$
\xi=g_*(t_1^0V(p)+ t_2^0W(p))+se_3(p),\;\;
\eta=g_*(t_1^0Y(p)+ t_2^0Z(p))+se_4(p).
$$
In particular, if $g$ is $1$-isotropic then 
$$
\|X_1\|=\Omega=\|X_2\|
\;\;\mbox{with}\;\;\<X_1,X_2\>=0\;\;\mbox{and}\;\;
\|\xi\|=\Omega=\|\eta\|\;\;\mbox{with}\;\;\<\xi,\eta\>=0
$$ 
where
$\Omega^2=s^2+\|t_1^0V(p)+t_2^0W(p)\|^2$.
\end{lemma}

\proof On one hand, 
$$
G_* Y_i=sg_*e_i(p)
+\sum_{j\geq3}\<\nap_{e_i}\delta_v,e_j\>(p)e_j(p),\;\;1\leq i\leq 2,
$$
gives
$$
G_*Y_i-\sum_{j\geq5}\<\nap_{e_i}\delta_v,e_j\>(p)G_*E_{j-2}
=sg_*e_i(p)-\sum_{3\leq k\leq 4}\<\nap_{e_i}e_k,\delta_v\>(p)e_k(p).
$$
On the other hand, 
\bea
\<\nap_{e_i}\delta_v,e_5\>(p)\!\!\!&=&\!\!\!-t_2^0\omega_{56}^i(p)
-t_3^0\omega_{57}^i(p)-t_4^0\omega_{58}^i(p)=-G_i(p,v),\\
\<\nap_{e_i}\delta_v, e_6\>(p)\!\!\!&=&\!\!\!t_1^0\omega_{56}^i(p)
-t_3^0\omega_{67}^i(p)-t_4^0\omega_{68}^i(p)=-H_i(p,v),\\
\<\nap_{e_i} e_3,\delta_v\>(p)\!\!\!&=&\!\!\!t_1^0\omega_{35}^i(p)+t_2^0\omega_{36}^i(p)
= t_1^0a_i(p)+t_2^0b_i(p),\\
\<\nap_{e_i} e_4,\delta_v\>(p)\!\!\!&=&\!\!\!t_1^0\omega_{45}^i(p)+t_2^0\omega_{46}^i(p)
= t_1^0c_i(p)+t_2^0d_i(p).
\eea
Hence, 
$$
G_*X_i=sg_* e_i-(t_1^0a_i+t_2^0b_i)e_3-(t_1^0c_i+t_2^0d_i)e_4,\;\; 1\leq i\leq 2.
$$
The remaining of the proof is straightforward using (\ref{omegas}).
\qed

\begin{lemma}\po\label{comp}
The following equations hold:
\be\label{zero}
\xi_*\d/\d s=e_3,\;\;\eta_*\d/\d s=e_4, 
\ee
\be\label{first}
\xi_* E_3=g_*V,\;\;\xi_* E_4=g_*W\;\;\mbox{and}\;\;\xi_*=0\;
\text{on}\;\mathcal{V}^0,
\ee
\be\label{second}
\eta_*E_3=g_*Y,\;\;\eta_* E_4=g_*Z\;\;\mbox{and}\;\;
\eta_*=0\;\text{on}\;\mathcal{V}^0,
\ee
\begin{eqnarray}\label{1}
\xi_* X_1\!\!\!&=&\!\!\!g_*\big((e_1(\va_1)-s\kappa)e_1+e_1(\va_2)e_2
+\omega_{12}^1J(t_1^0V+t_2^0W)+G_1V+H_1W\big)\nonumber\\
\!\!\!&&\!\!\!+\kappa\va_1e_3+(s\omega_{34}^1
+\lambda\kappa\va_2)e_4+sa_1e_5+sb_1e_6-\va_1 g,
\end{eqnarray}
\begin{eqnarray}\label{2}
\xi_* X_2\!\!\!&=&\!\!\!g_*\big(e_2(\va_1)e_1+(e_2(\va_2)+s\kappa)e_2
+\omega_{12}^2J(t_1^0V+t_2^0W)+G_2V+H_2W\big)\nonumber\\
\!\!\!&&\!\!\!-\kappa\va_2e_3+(s\omega_{34}^2
+\lambda\kappa\va_1)e_4+sa_2e_5+sb_2e_6-\va_2 g,
\end{eqnarray}
\begin{eqnarray}\label{3}
\eta_* X_1\!\!\!&=&\!\!\!g_*\big(e_1(\psi_1)e_1+(e_1(\psi_2)-s\lambda\kappa)e_2
+\sigma\omega_{12}^1(t_1^0V+t_2^0W)-\sigma G_1JV-\sigma H_1JW\big)\nonumber\\
\!\!\!&&\!\!\!-(s\omega_{34}^1-\kappa\psi_1)e_3
+\lambda\kappa\psi_2e_4+s\sigma a_2e_5+s\sigma b_2e_6-\psi_1 g,
\end{eqnarray}
\begin{eqnarray}\label{4}
\eta_* X_2\!\!\!&=&\!\!\!g_*\big((e_2(\psi_1)-s\lambda\kappa)e_1+e_2(\psi_2)e_2
+\sigma\omega_{12}^2(t_1^0V+t_2^0W)-\sigma G_2JV-\sigma H_2JW
\big)\nonumber\\
\!\!\!&&\!\!\!-(s\omega_{34}^2+\kappa\psi_2)e_3
+\lambda\kappa\psi_1e_4-s\sigma a_1e_5-s\sigma b_1e_6-\psi_2 g
\end{eqnarray}
where $\sigma=1/\lambda$ and 
$\psi_j=t_1^0c_j+t_2^0d_j,\; j=1,2$.
\end{lemma}

\proof We compute at $(s,p,v)\in N^{n+1}$.
Let $\gamma(t)=(s,p,v(t))$ be a curve in $N^{n+1}$ such that 
$v(0)=v$, and thus $\gamma'(0)\in\mathcal{V}(p,v)$. We have that
$$
\xi_*\gamma'(0)=\<Dv/dt(0),e_5(p)\>g_*V(p)+\<Dv/dt(0),e_6(p)\> g_*W(p), 
$$
or equivalently, that
$$
\xi_*\gamma'(0)=\<G_*\gamma'(0),e_5(p)\>g_*V(p)
+\<G_*\gamma'(0),e_6(p)\> g_*W(p). 
$$
 From this we obtain (\ref{first}). Similarly, we have (\ref{second}). 

To obtain (\ref{1}) to (\ref{4}) one has to use  Lemma \ref{NF} and 
the Gauss and Weingarten  formulas for $g$. 
We only argue for  (\ref{1})  since the proof of the other  equations 
is similar.  We have from (\ref{X}) and (\ref{first}) that
$$
\xi_*X_i =\xi_*Y_i+G_ig_*V +H_ig_*W,\;\; 1\leq i\leq 2. 
$$
In view of (\ref{yes}) and since
$$
\xi\circ\beta_i(t)=t_1^0g_*V(c_i(t))
+t_2^0g_*W(c_i(t))+se_3(c_i(t)),
$$
we obtain in terms of the connection in $L^2$ that 
\bea
\xi_*Y_i\!\!\!&=&\!\!\!t_1^0\big(g_*\n_{e_i}V +\a_g(e_i,V)\big)(p)
+t_2^0\big(g_*\n_{e_i}W +\a_g(e_i,W)\big)(p)\\
\!\!\!&&\!\!\!+(-1)^is\kappa(p)g_*e_i(p)+s\nap_{e_i}e_3(p),
\eea
and (\ref{1}) follows by a direct computation.
\qed

\begin{lemma}\po\label{A} The second fundamental form of $G$ 
in terms of the orthonormal frame 
$$
E_0= \d/\d s,\;\;E_i=X_i/\Omega,\;i=1,2,\;\;\mbox{and}\;\;G_*E_j=e_{j+2},\;3\leq j\leq n,
$$
vanishes along $\mathcal{V}^0$ and restricted to 
$\mathrm{span}\{E_0\}\oplus\mathcal{H}^G\oplus\mathcal{V}^1$  
is given by
$$
A_{\xi}=\begin{bmatrix}
0&\bar\va_1&\bar\va_2&0&0\\
\bar\va_1&h_1+\kappa&h_2&r_1&s_1\\
\bar\va_2&h_2&-h_1-\kappa&r_2&s_2\\
0&r_1&r_2&0&0\\
0&s_1&s_2&0&0&
\!\!\!\!\!\end{bmatrix},\;
A_{\eta}=\begin{bmatrix}
0&\bar\va_2&-\bar\va_1&0&0\\
\bar\va_2&h_2&\kappa-h_1&r_2&s_2\\
-\bar\va_1&\kappa-h_1&-h_2&-r_1&-s_1\\
0&r_2&-r_1&0&0\\
0&s_2&-s_1&0&0&
\!\!\!\!\!\end{bmatrix}\\
$$
where $\bar\va_i\Omega=\va_i$, $r_i\Omega=-sa_i$, $s_i\Omega=-sb_i$ and 
\bea
h_i\!\!\!&=&\!\!\!-\dfrac{s}{\Omega^2}\big(t_1(e_i(a_1)
- a_2B_i- b_1\omega_{56}^i) 
+ t_2(e_i( b_1)-b_2B_i+ a_1\omega_{56}^i)\\
\!\!\!&&\!\!\! +\,t_3(a_1\omega_{57}^i+ b_1\omega_{67}^i)
+t_4(a_1\omega_{58}^i+ b_1\omega_{68}^i)\big)
\eea
with $B_i=\omega_{12}^i+\omega_{34}^i$,\; $i=1,2$.
\end{lemma}

\proof Since $g$ is $1$-isotropic, then (\ref{1}) to (\ref{4}) 
hold for $\psi_1=\va_2$ and $\psi_2=-\va_1$.
On the other hand, a straightforward computation shows that 
the Ricci equations  
$$
\<R^\perp(e_1,e_2) e_\a,e_\beta\>=0
$$
for $\a=3,4$ and $\beta=5,6$
are equivalent to
\bea
&&e_1(a_2)-e_2(a_1)+a_1B_1+a_2B_2-b_2\omega_{56}^1 
+b_1\omega_{56}^2=0,\\
&&e_1(b_2)-e_2(b_1)+b_1B_1+b_2B_2+a_2\omega_{56}^1
-a_1\omega_{56}^2=0,\\
&&e_1(a_1)+e_2(a_2)-a_2B_1+a_1B_2-b_1\omega_{56}^1
-b_2\omega_{56}^2=0,\\
&&e_1(b_1)+e_2(b_2)-b_2B_1+b_1B_2+a_1\omega_{56}^1
+a_2\omega_{56}^2=0,
\eea
and for $\a=3,4$ and $\beta=7,8$ are equivalent to
\bea
&&a_2\omega_{57}^1-a_1\omega_{57}^2 
+b_2\omega_{67}^1-b_1\omega_{67}^2=0,\\
&&a_2\omega_{58}^1-a_1\omega_{58}^2 
+b_2\omega_{68}^1-b_1\omega_{68}^2=0,\\
&&a_1\omega_{57}^1+a_2\omega_{57}^2 
+b_1\omega_{67}^1+b_2\omega_{67}^2=0,\\
&&a_1\omega_{58}^1+a_2\omega_{58}^2 
+b_1\omega_{68}^1+b_2\omega_{68}^2=0.
\eea
We thus have that
$$
\<A_\xi E_i,E_j\>=-\<G_*E_i,\xi_* E_j \>\;\;\mbox{and}\;\; 
\< A_\eta E_i,{E_j}\>=-\< G_*E_i ,\eta_* E_j \>,\;0\leq i,j \leq n,
$$ 
and the result follows by a straightforward computation.\vspace{1,5ex}\qed

\noindent\emph{Proof of Theorem \ref{main1}:}
We first prove the converse. If $F\colon M^n\to\Sf^{n+2}$, $n\geq 4$, 
is an $(n-2)$-ruled minimal immersion with rank $\rho=4$ everywhere, 
then the tangent bundle splits as $TM={\cal H}\oplus{\cal V}$, 
where $\cal H$ is orthogonal to the rulings. Moreover,  we have that ${\cal V}$ 
splits as ${\cal V}={\cal V}^1\oplus{\cal V}^0$ with the fibers of ${\cal V}^0$ 
being the relative nullity leaves.  

 The normal space of the surface  $g=F\circ j$ at $x\in L^2$ is given by 
$$
N_gL(x)=F_*(j(x)){\cal V}\oplus N_FM(j(x)).
$$
Being $j$ is totally geodesic, we have 
\be\label{ag}
\a_g(X,Y)=\a_F(j_*X,j_*Y)
\ee
for all $X,Y\in TL$. This and our assumptions imply that $g$ is minimal. 

Let $\pi\colon\Lambda_g\to L^2$ denote the subbundle of the normal bundle of $g$ 
whose fiber at $x\in L^2$ is $F_*(j(x)){\cal V}$.   
We consider the cone $\mathcal C F\colon\R\times M^n\to\R^{n+3}$ 
given by 
$$
\mathcal CF(t,p)=tF(p).
$$ 
Observe that 
$$
\mathcal C F (t,p)-\mathcal C F (u(t,p),j(x))=\mathcal C F (t,p)
-u(t,p)g\circ\pi(p)\in F_*(j(x)){\cal V}
$$
for any $p\in M^n$, where $x=\pi(p)$, since $p$ and $j(x)$ belong to the 
same leaf of  ${\cal V}$ and
$$
u(t,p)=t/\<F(p),g\circ\pi (p)\>.
$$
Since $\mathcal C F$ maps locally diffeomorphically the leaves of ${\cal V}$  
onto   affine subspaces, it follows that the map 
$T \colon \R \times M^n \to\R \times\Lambda_g$ given by
$$
T(t,p)=(u(t,p), \pi(p),\mathcal C F (t,p)-u(t,p) g\circ \pi(p)) 
$$
is  a local diffeomorphism.  Clearly the  immersion 
$\tilde{G}=\mathcal C F\circ T^{-1}$ satisfies 
$$
\tilde{G}(s,x,v)=sg(x)+v,
$$
i.e., $\tilde{G}=G_g$ is of the form  (\ref{map}). 
Identifying locally $ \R \times M^n$ with $ \R \times\Lambda_g$ via $T$, 
we have that $\mathcal C F=G_g=G$ and $j$ is the zero section of $\Lambda_g$, 
i.e., we have the parametrization given by (\ref{map}). 
The horizontal and the vertical bundles satisfy
$$
G_*(s,p,v)\mathcal V=(N_1^g(p))^\perp\subset N_gL(p),\;\;
G_*(s,p,v)\mathcal H^G\subset g_*T_pL\oplus(\Lambda_g(p))^\perp,
$$
$$
N_GN(s,p,v)\subset g_*T_pL\oplus(\Lambda_g(p))^\perp
$$
and now (\ref{ag}) yields  $N_1^g  = \Lambda_g^\perp$.   

It remains to see that $g$ is $1$-isotropic.
For an adapted frame  $\{e_1,\ldots,e_{n+2}\}$ of $g$ set 
$$
g_{ij}=\<G_*X_i,G_*X_j\>
$$
and 
$$
b^\xi_{ij}=\<\xi_*X_i,G_*X_j\>,\;\;\; 
b^\eta_{ij}=\<\eta_*X_i,G_*X_j\>,\;\;i,j=1,2.
$$
Using Lemma \ref{NF} and Lemma \ref{comp}, we find that
$$
g_{11}=s^2+\va_1^2+\sigma^2\va_2^2,\;\;
g_{12}=(1-\sigma^2)\va_1\va_2,\;\;
g_{22}=s^2+\va_2^2+\sigma^2\va_1^2,
$$
and
\bea
b^\xi_{11}\!\!\!&=&\!\!\!s(e_1(\va_1)-s\kappa-\omega_{12}^1\va_2+G_1a_1+H_1b_1)
-\kappa\va_1^2-\sigma\va_2(s\omega_{34}^1+\mu\va_2),\\
b^\xi_{12}\!\!\!&=&\!\!\!s(e_1(\va_2)+\omega_{12}^1\va_1+G_1a_2+H_1b_2)
-\kappa\va_1\va_2 +\sigma\va_1(s\omega_{34}^1+\mu\va_2),\\
b^\xi_{21}\!\!\!&=&\!\!\!s(e_2(\va_1)-\omega_{12}^2\va_2+G_2a_1+H_2b_1)
+\kappa\va_1\va_2-\sigma\va_2(s\omega_{34}^2+\mu\va_1),\\
b^\xi_{22}\!\!\!&=&\!\!\!s(e_2(\va_2)+s\kappa+\omega_{12}^2\va_1+G_2a_2+H_2b_2)
+\kappa\va_2^2+\sigma\va_1(s\omega_{34}^2+\mu\va_1)
\eea
and
\bea
b^\eta_{11}\!\!\!&=&\!\!\!s(e_1(\psi_1)-\omega_{12}^1\psi_2+\sigma G_1a_2+\sigma H_1 b_2)
+s\omega_{34}^1\va_1-\kappa(\va_1\psi_1+\va_2\psi_2),\\
b^\eta_{12}\!\!\!&=&\!\!\!s(e_1(\psi_2)-\mu+\omega_{12}^1\psi_1- \sigma G_1a_1-\sigma H_1b_1)
+s\omega_{34}^1\va_2+\kappa(\va_1\psi_2-\va_2\psi_1),\\
b^\eta_{21}\!\!\!&=&\!\!\!s(e_2(\psi_1)-\mu-\omega_{12}^2\psi_2+\sigma G_2a_2+\sigma H_2b_2)
+s\omega_{34}^2\va_1+\kappa(\va_1 \psi_2-\va_2\psi_1),\\
b^\eta_{22}\!\!\!&=&\!\!\!s(e_2(\psi_2)+\omega_{12}^2\psi_1-\sigma G_2a_1-\sigma H_2b_1)
+s\omega_{34}^2\va_2 +\kappa(\va_1\psi_1+\va_2\psi_2).
\eea
 From our assumptions, we have 
$$
g_{11}b^\xi_{22}-g_{12}(b^\xi_{12}+b^\xi_{21})+g_{22}b^\xi_{11}=0\;\;\mbox{and}\;\;
g_{11}b^\eta_{22}-g_{12}(b^\eta_{12}+b^\eta_{21})+g_{22}b^\eta_{11}=0.
$$
Viewing these as polynomials were the coefficients of 
$t_1^4, t_2^4$ and $t_1^2t_2^2$  must vanish gives
$$
(\lambda^2-1)(a_1^2+a_2^2)(a_1^2-a_2^2)=0=(\lambda^2-1)(b_1^2+b_2^2)(b_1^2-b_2^2)
$$
and
$$
(\lambda^2-1)a_1a_2(a_1^2+a_2^2)=0=(\lambda^2-1)b_1b_2(b_1^2+b_2^2).
$$
Hence $\lambda=1$ since, otherwise, we would have that
$\omega_{35}=\omega_{36}=\omega_{45}=\omega_{46}=0$, 
and that is a contradiction.\vspace{1ex} 

We now prove the direct statement. Since $F_g=G|_M$, we obtain that
$g=F_g\circ j$ where $j\colon L^2\to M^n$ is given by $j(x)=(\pm1,x,0)$.  
Clearly, we have that $j$  is an  integral surface of the distribution 
orthogonal to the rulings that is  totally geodesic and a global 
cross section to the rulings.  Up to uniqueness of the integral 
surface  the  proof follows from Lemma \ref{A}. 

Assume  that there exists a second integral surface 
$\tilde j\colon L^2\to M^n$. 
Set $\tilde g=F_g\circ\tilde j$ and let  
$\tilde T\colon\R\times M^n\to\R\times\Lambda_{\tilde g}$ 
be the local diffeomorphism given by 
$$
\tilde T(t,p)=(\tilde u (t,p), \pi(p),\mathcal CF(p)
-\tilde u (t,p)\tilde g\circ \pi(p))
$$
where
$$
\tilde u(t,p)=t/\<F(p),\tilde g\circ\pi (p)\>.
$$
Then  $\tilde T\circ T^{-1}\colon \R \times\Lambda_g\to\R\times\Lambda_{\tilde g}$ 
is given by
$$
\tilde T\circ T^{-1} (s,x,v)=(\tilde s, x,v+sg(x)- \tilde s \tilde g(x)),
$$
where $T^{-1} (s,x,v)=(t,p)$ and $\tilde s = \tilde u (t,p)$. 
Hence $\Lambda_g$ and $\Lambda_{\tilde g}$  can be identified by parallel 
translation. Using that $sg(x)- \tilde s \tilde g(x) \in \Lambda_g(x)$, 
we obtain that $\tilde g = \pm g$.
\vspace{1,5ex}\qed

The vertical bundle $\mathcal V$ of the submersion $\pi$ given by 
$\mathcal V=\ker\pi_*$  can be orthogonally decomposed as
$\mathcal{V}=\mathcal{V}^1\oplus \mathcal{V}^0$ on an open dense subset 
of $L^2$, where $\mathcal{V}^1$ denotes the plane bundle determined by $N_2^g$.
In fact, this holds if $N_1^g$ and $N_2^g$ are subbundles, which we can assume 
without loss of generality.
In the sequel, we consider the orthogonal decomposition of the tangent bundle
of $N^{n+1}$ given by 
$$
TN=\mathrm{span}\{\d/\d s\}\oplus{\cal H}^G\oplus\mathcal{V}
$$
where we identify isometrically (and use the same
notation) the subbundle $\mathcal{V}$ tangent to the rulings with the 
corresponding normal subbundle to $g$. Then, it follows from the proof
that the relative nullity leaves of $G$ are identified with the fibers 
of $\mathrm{span}\{\d/\d s\}\oplus\mathcal{V}^0$.

\vspace{1,5ex}

Let $\cal J$ denote the endomorphism such that 
${{\cal J}|_{{\cal H}^G}}\colon{\cal H}^G\to{\cal H}^G$ is the almost complex 
structure in $\mathcal{H}^G$ determined by the orientation and restricted 
to $\mathrm{span}\{\d/\d s\}\oplus\cal V$ is the identity and set
$$
{\cal J}_{\theta}=\cos\theta I+\sin\theta {\cal J}.
$$

\noindent\emph{Proof of Theorem \ref{af}:} 
For each $\theta\in \Sf^1$ consider the
submanifold $G_\theta \colon N^{n+1}\to\R^{n+3}$ defined by 
$$
G_\theta(s,p,v)=sg_\theta (p)+\phi_\theta v
$$
where $\phi_\theta\colon N_gL\to N_{g_{\theta}}L$ is the parallel 
vector bundle isometry that identifies the normal subbundles of $g$ 
and of $g_\theta$.

 In the sequel, corresponding quantities of $G_\theta$ are denoted 
by the same symbol used for $G$ marked with $\theta$.  
That  $G_\theta$ is isometric to $G$ is immediate.
Since the tangent frame $\{e_1,e_2\}$
has been fixed, we have for the adapted frames of $g_\theta$ that
$$
e^\theta_3=\phi_\theta\circ R^1_\theta e_3\;\;\mbox{and}\;\; 
e^\theta_4=\phi_\theta\circ R^1_\theta e_4        
$$
where $R^1_\theta$ is the  rotation of angle $\theta$ on $N^g_1$. 
We complete the adapted frame choosing 
$$
e^\theta_j=\phi_\theta e_j,\;\;\; 5\leq j \leq n+2.
$$
Clearly, it holds that $\omega^\theta_{34}=\omega_{34}$ 
and $\omega^\theta_{ij}=\omega_{ij}$ for $i,j\geq 5$. Moreover, 
$$
\omega^\theta_{35}
=\cos\theta\omega_{35}-\sin\theta*\omega_{35}\;\;\mbox{and}\;\;
\omega^\theta_{36}=\cos\theta\omega_{36}-\sin\theta*\omega_{36}.
$$
Hence, the dual vector fields of $\omega^\theta_{36}$  
and $\omega^\theta_{36}$ are given, respectively, by
$$
V_\theta =J_{-\theta}V \;\; \text{and}\;\; W_\theta =J_{-\theta}W.
$$
Thus,
$$
a_1^\theta= a_1\cos\theta+  a_2\sin\theta,\;\;a_2^\theta
= a_2\cos\theta- a_1\sin\theta,
$$
$$
b_1^\theta= b_1\cos\theta+b_2\sin\theta,\;\;b_2^\theta
=b_2\cos\theta-b_1\sin\theta.
$$
It follows from (\ref{delta}), (\ref{yes}) and (\ref{X}) that 
$$
X_i^\theta=X_i,\;\;i=1,2.
$$
By Lemma \ref{NF}, the normal bundle of $G_\theta$ is spanned by 
$$
\xi_\theta=g_{\theta_*}J_{-\theta}(t_1V+ t_2W)+s\phi_\theta\circ R^1_\theta e_3,\;\;\;
\eta_\theta=-g_{\theta_*}J_{\pi/2-\theta}(t_1V+t_2W)+s\phi_\theta\circ R^1_\theta e_4.
$$
A straightforward computation yields that the map
$\Psi_\theta\colon N_GN\to N_{G_\theta}N$ given by
$$
\Psi_\theta \xi=\xi_\theta \;\;\text{and}\;\;\Psi_\theta\eta=\eta_\theta
$$
is a parallel vector bundle isometry.  The shape operators of 
$G_\theta$ vanish on $\mathcal{V}^0$
and restricted to $\mathrm{span}\{\d/\d s\}\oplus\mathcal{H}^G\oplus\mathcal{V}^1$ 
and with respect to $\{E_1,\dots,E_n\}$ they are given by
$$
A^\theta_{\xi_\theta}=\begin{bmatrix}
0&\bar\va_1^\theta&\bar\va_2^\theta&0&0\\
\bar\va_1^\theta&h_1^\theta+\kappa&h_2^\theta&r_1^\theta&s_1^\theta\\
\bar\va_2^\theta&h_2^\theta&-h_1^\theta-\kappa&r_2^\theta&s_2^\theta\\
0&r_1^\theta&r_2^\theta&0&0\\
0&s_1^\theta&s_2^\theta&0&0&
\!\!\!\!\!\end{bmatrix}, 
\;
A^\theta_{\eta_\theta}=\begin{bmatrix}
0&\bar\va_2^\theta&-\bar\va_1^\theta&0&0\\
\bar\va_2^\theta&h_2^\theta&\kappa-h_1^\theta&r_2^\theta&s_2^\theta\\
-\bar\va_1^\theta&\kappa-h_1^\theta&-h_2^\theta&-r_1^\theta&-s_1^\theta\\
0&r_2^\theta&-r_1^\theta&0&0\\
0&s_2^\theta&-s_1^\theta&0&0&
\!\!\!\!\!\end{bmatrix}
$$
where $\bar\va_i^\theta\Omega=\va_i^\theta$, $r_i^\theta\Omega=-sa_i^\theta$, 
$s_i^\theta\Omega=-sb_i^\theta$\; $i=1,2$,  and 
$$
\va_1^\theta=\va_1 \cos\theta +\va_2 \sin\theta,
\;\;\;\va_2^\theta= -\va_1\sin\theta+\va_2\cos\theta,
$$
$$
h_1^\theta=h_1 \cos\theta +h_2 \sin\theta,
\;\;\;h_2^\theta= -h_1\sin\theta+h_2\cos\theta.
$$

Let $L_\theta\colon TN\to TN$ be such that 
$L_\theta |_{\mathrm{span}\{\d/\d s\}\oplus\cal{V}}=0$ and
$L_\theta|_{{\cal H}^G}\colon{\cal H}^G\to{\cal H}^G$ is 
the reflection given by  
$$
L_\theta|_{{\cal H} ^G}=\begin{bmatrix}
-\sin(\theta/2)&\cos(\theta/2)\\
\cos(\theta/2)&\sin(\theta/2)&
\!\!\!\!\!\end{bmatrix}
$$
with respect to the tangent frame $\{E_1,E_2\}$.  It follows easily  that
$$
A^\theta_{\Psi_\theta\xi}=A_{\textsf{R}_\theta\xi} 
-2\kappa\sin(\theta/2)L_\theta\;\;\mbox{and}
\;\;A^\theta_{\Psi_\theta\eta}=A_{\textsf{R}_\theta\eta}
-2\kappa\sin(\theta/2){\cal J}\circ L_\theta.
$$
By direct computation, we obtain
$$
\a_{G_\theta}(X,Y)=\Psi_\theta \Big(\textsf{R}_{-\theta} \a_G(X,Y) 
-\frac{2\kappa}{\Omega^2}\sin(\theta/2) (\<L_\theta X,Y\>\xi
+\<L_\theta {\cal J}X,Y\>\eta\Big).
$$
Now let  $\beta$ be the symmetric section of 
$Hom(TN\times TN,N_GN)$ with nullity $\mathcal V$ given by
\be\label{b}
\beta(E_1,E_1)=\frac{1}{\Omega^2}\xi=-\beta(E_2,E_2),
\;\;\beta(E_1,E_2)=-\frac{1}{\Omega^2}\eta,
\ee
and the proof of (\ref{forms}) follows easily.

Finally, that the isometric deformations $G_\theta$ of $G$ are genuine 
is immediate from Lemma \ref{A} since the shape operators of $G$ 
have rank four for any normal direction along an open dense subset of $N^{n+1}$. 
\vspace{1,5ex}\qed

\noindent\emph{Proof of Theorem \ref{main2}:} Given $\theta\in\Sf^1$, denote 
$F_\theta=G_\theta|_M$ where
$$
G_\theta(s,p,v)=sg_\theta(p)+\phi_\theta v.
$$
That $F_g$ allows a one-parameter family of minimal isometric immersions 
$F_\theta\colon M^n\to\Sf^{n+2}$, $\theta\in\Sf^1$, such that $F_0=F_g$ 
and each $F_\theta$ carries the same ruling and relative nullity leaves 
as $F_g$ is a consequence of Proposition \ref{af}.\qed
\vspace{1,5ex}

\noindent\emph{Proof of Theorem \ref{main3}:} 
It is completely analogous to the proof of Theorem $6$ in \cite{dv1}.\qed
\vspace{1,5ex}

\noindent\emph{Proof of Theorem \ref{mod}:}
Let $\bar{F}\colon M^3\to\Sf^5$ be a ruled isometric minimal immersion 
with the same rulings  as $F_g$ and set $\bar g=\bar F\circ j$.  From the 
proof of Theorem \ref{main1}, we have that the surface 
$\bar{g}$ is isometric to $g$ and isotropic.  Hence, the set of all minimal 
isometric immersions of $M^3$ into $\Sf^5$ with the same rulings as $F_g$ 
can be identified with the set of all  isotropic immersions of $L^2$ into
$\Sf^5$. The proof now follows from the results in \cite{dv2}.\qed
\vspace{1,5ex}

\noindent\emph{Proof of Theorem \ref{last}:} 
Using Lemma \ref{NF} and Lemma \ref{A}, we have that the
squared length of the second fundamental form of $G$ is given by
$$
\|\a_G\|^2=\dfrac{4}{\Omega^4}\Big((1-K)\Omega^2+\va_1^2
+\va_2^2+\Omega^2\sum_{i=1,2}(h_i^2+r_i^2+s_i^2)\Big).
$$
It follows that
\be\label{length}
\|\a_G\|^2(s,p,v)=\dfrac{4}{\Omega^2}\Big(2-K+h_1^2+h_2^2
+\dfrac{s^2}{\Omega^2}(\|V\|^2+\|W\|^2-1)\Big).
\ee

By Corollary $4$ in \cite{V} any $1$-isotropic  torus in $\Sf^5$ 
is regular, hence $M^3$ is compact.
On the other hand, we have that $g$ is $O(6)$-congruent to a holomorphic 
curve in the nearly Kaehler  sphere  $\Sf^6$; see \cite{EV} or \cite{V}. 
Choose local orthonormal frame 
$\{e_1,e_2,e_3,e_4, e_5\}$  
such that 
$$
\a_g(e_1,e_1)=\sqrt{1/2} e_3,\;\;\; \a_g(e_1,e_2)=\sqrt{1/2} e_4,
$$ 
$$
\a^3_g(e_1,e_1,e_1)=\kappa_1 e_5,\;\;\; \a^3_g(e_1,e_1,e_2)=0,
$$
where $\kappa_1=\sqrt{1/2}$ by Theorem 5 in \cite{V}. 
Hence, we have that $V=e_1$. From Lemma $6$ in \cite{V0} 
we obtain $h_1=h_2=0$.  Now (\ref{length}) gives
$$
\|\a_G\|^2(s,p,v)=\dfrac{8}{\Omega^2}=\dfrac{8}{s^2+t_1^2},
$$
and hence $\|\a_F\|^2=8$.\qed

{\renewcommand{\baselinestretch}{1}
\hspace*{-20ex}\begin{tabbing} \indent\= IMPA -- Estrada Dona Castorina, 110
\indent\indent\= Univ. of Ioannina -- Math. Dept. \\
\> 22460-320 -- Rio de Janeiro -- Brazil  \>
45110 Ioannina -- Greece \\
\> E-mail: marcos@impa.br \> E-mail: tvlachos@uoi.gr
\end{tabbing}}
\end{document}